\documentclass[12pt]{article}%
\usepackage{amsmath}
\usepackage{amsmath}
\usepackage{amssymb}
\usepackage{dsfont}
\usepackage{cite}
\newsymbol \blackbox 1004
\newcommand{\eh}{\hfill}\newlength{\sperr}

\newenvironment{proof}{{\settowidth{\sperr}{\bf\rm
Proof}%
\par\addvspace{0.3cm}\noindent\parbox[t]{1.3\sperr}
{\bf\rm P\eh r\eh o\eh o\eh f\eh }%
}}{\nopagebreak\mbox{}
$\blackbox$\par\addvspace{0.3cm}}

\def\nn{\nonumber}

\def\b{\beta}

\def\s{\sigma}
\def\la{\lambda}

\def\ze{\zeta}

\def\vt{\vartheta}

\def\wt{\widetilde}

\def\p{\partial}

\def\BC{{\mathbb C}}

\def\BR{{\mathbb R}}
\def\BN{{\mathbb N}}

\def\cla{{\mathcal A}}
\def\clb{{\mathcal B}}

\def\clj{{\mathcal J}}

\def\cli{{\mathcal I}}

\def\cln{{\mathcal N}}

\def\cls{\mathcal{S}}
\def\clt{\mathcal{T}}

\def\diag{\mathrm{diag}}
\newcommand{\E}{\mathrm{e}}
\newcommand{\I}{\mathrm{i}}

\newtheorem{Pa}{Paper}[section]
\newtheorem{Tm}[Pa]{{\bf Theorem}}

\newtheorem{Cy}[Pa]{{\bf Corollary}}
\newtheorem{Rk}[Pa]{{\bf Remark}}

\newtheorem{Ee}[Pa]{{\bf Example}}

\newtheorem{Pn}[Pa]{{\bf Proposition}}

\title{Dressing for generalised linear  Hamiltonian systems
depending rationally on the spectral parameter and some applications}

\author{Alexander Sakhnovich}

\date{}
\begin{document}
\maketitle

\begin{flushright}
Faculty of Mathematics,
University
of
Vienna, \\
Oskar-Morgenstern-Platz 1, A-1090 Vienna,
Austria, \\
e-mail: {\tt oleksandr.sakhnovych@univie.ac.at}
\end{flushright}

\begin{abstract}  
We construct so called Darboux matrices and fundamental solutions
in the important case of  the generalised Hamiltonian (or canonical) systems
depending rationally on the spectral parameter. A wide
class of explicit solutions is obtained in this way.
Interesting results for dynamical systems depending on several variables
and their explicit solutions follow.
For these purposes we use
our version of B\"acklund-Darboux transformation
and square roots of  the corresponding generalised matrix eigenvalues.
Some new auxiliary  results on the roots of  matrices are 
included as well.  
An appendix is added to make the paper self-sufficient.
\end{abstract}


\vspace{0.2em}


\section{Introduction} \label{intro}
\setcounter{equation}{0}
{\bf 1.} In this paper, we construct fundamental solutions
and so called Darboux matrices for generalised Hamiltonian (canonical) systems
depending rationally on the spectral parameter:
\begin{align} &       \label{I1}
y^{\prime}(x,z)=\I \sum_{k=1}^{r} (z-c_k)^{-1}j\b_k(x)^*\b_k(x)y(x,z) \quad \left(y(x,z)\in \BC^m\right);
\\ &       \label{I2}
   j=\begin{bmatrix} I_{m_1} & 0\\ 0 & -I_{m_2}\end{bmatrix},  \,\, \b_k(x)\in \BC^{p\times m} \quad (p\leq m), \,\, c_k\in \BR; \,\,  m_1,m_2\in \BN_0,
\\ &       \label{I2'} 
\BN_0=\BN\cup \{0\}, \quad m_1+m_2=m>0.  
\end{align}
Here, $\I$ is the  imaginary unit ($\I^2=-1$), the $m\times m$ Hamiltonians
$H_k(x)=\b_k(x)^*\b_k(x)$ are locally summable, $c_i\not=c_k$ for $i\not=k$, $\BR$ stands for the real axis,  $\BC$ stands for the complex plane,
$\BN$ stands for the set of positive integers,
$I_n$ denotes the $n\times n$ identity matrix, the class of $p\times m$ matrices with complex-valued  entries
is denoted by $\BC^{p\times m}$, and $\BC^{p}=\BC^{p\times 1}$.  We note that $j=I_m$ if $m_2=0$ and $j=-I_m$ if $m_1=0$.

The obtained results on systems \eqref{I1} are used in order to
construct solutions
of the dynamical systems
\begin{align} &       \label{I3}
\frac{\p \psi}{\p x}=\I\sum_{k=1}^r j H_k(x)\frac{\p \psi}{\p \ze_k}\,.
\end{align}
When $r=1$ and $H_1(x)$ is invertible, system \eqref{I3} is equivalent to the important dynamical canonical system 
studied, for instance,  in \cite{JacZ, ALS17}. Thus, this paper is a substantial development and generalisation of
the work \cite{ALS17}. The so called transformed Hamiltonians $\wt H_k$ treated here
preserve the property which the Hamiltonian in \cite{ALS17} possessed. Namely, the matrices $j\wt H_k(x)$ are linear similar to $jH_k(x)$. 
 (This property is essential for the research on the
$C_0$-semigroups related to dynamical canonical systems \cite{JacZ, ALS17}.)

Symmetries and conservation laws are important in the study of  the dynamical systems (see, e.g., \cite{Col, Co, Daf, HeHo}
and references therein). The expression for energy from \cite[Proposition 3.1]{ALS17} is not valid for system \eqref{I3}
but an interesting conservation law is presented in Proposition \ref{PnCL}.

{\bf 2.}  In order to construct fundamental solutions
and Darboux matrices for systems \eqref{I1},
we use
our GBDT (generalised B\"acklund-Darboux transformation) approach
and square roots of  {\it generalised matrix eigenvalues}.  The case $r=2$, $p=1$ and $j=I_2$ was studied earlier
in \cite{MST}. Here, we greatly generalise the corresponding results from \cite{MST}.

Dressing procedure (also called B\"acklund-Darboux transformation or commutation method)  is a fruitful tool in the study of linear and nonlinear differential equations
(see  \cite{Cies3, Deift, Ge, GeHo, GeT,  Gu,  Mar, MS, Mi, ZM} and numerous references therein).
For the recent applications and references see, for instance, \cite{CoI, Cu, Gar, Grun, ALS17, ALSJDE21}.
Dressing procedure may be considered as an important part of the general symmetry theory
\cite{Olv, RP}.
GBDT version of B\"acklund-Darboux transformation (with {generalised matrix eigenvalues}) was actively developed starting 
from our work \cite{SaA94} (see further references in \cite{Cies3, FKS, KaS, KoSaTe,   ALSa1, ALSstring, ALSJDE21, SaSaR}).

The notion of Darboux matrix (matrix function) may be explained in the following way. Given some initial linear system
$y^{\prime}(x,z)=G(x,z)y(x,z)$ $\left(y^{\prime}:=\frac{d}{dx}y\right)$ and a transformed (via B\"acklund-Darboux transformation) 
system $\wt y^{\prime}(x,z)=\wt G(x,z)\wt y(x,z)$,
the corresponding Darboux matrix is a nondegenerate matrix function  $W(x,z)$ satisfying equation 
\begin{align} &       \label{I0}
W^{\prime}(x,z)=\wt G(x,z) W(x,z)-W(x,z)G(x,z).
\end{align}
It is easy to see that the product $\wt w(x,z)=W(x,z)w(x,z)$, where $w(x,z)$ is a fundamental solution  of the initial system $y^{\prime}(x,z)=G(x,z)y(x,z)$,
is a fundamental solution  of the transformed system
$\wt y^{\prime}(x,z)=\wt G(x,z)\wt y(x,z)$.
If $G(x,z)$ has some simple form so that $w(x,z)$ and Darboux matrix $W(x,z)$ are known explicitly, then the fundamental solution  $\wt w(x,z)$ of the transformed
system
 is also known explicitly.

{\bf 3. } In Section \ref{gencan}, we construct GBDT for system \eqref{I1}, \eqref{I2}.
If  $r=1$ and $m_1=m_2$, system  \eqref{I1} is equivalent to the
classical canonical system (see various references and results on canonical systems in \cite{dB1, GoKr,  SaSaR, SaL2}).
L.~de~Branges and M.G. Krein obtained some fundamental results for the case of canonical systems where $m=2$
(see, e.g., \cite{dB1, GoKr}). The spectral theory of  canonical systems, where $m>2$, is no less important but
essentially more complicated (see, e.g., \cite{Ach, dB2, ALSstring, SaSaR, SaL2, Sep} and references
therein). In Section~\ref{gencan} we deal with the case of an arbitrary $m\in \BN$.

In the classical works and classical problems, there is a linear dependence on the spectral parameter
(see also some recent papers, for instance, \cite{Ach, GeS, Zem}).
The case of the more complicated (including nonlinear) dependence on the spectral parameter $z$ is of essential interest and was studied, for instance,
in \cite{Ful, ElH, HMS, MST,  MScS, Zhou}. Various other examples include auxiliary linear systems for  integrable nonlinear
equations. (See, e.g.,  \cite{MST} on the connection of the system \eqref{I1}, where $r=m=2$ and $p=1$, with the sine-Gordon equation.)

{\bf 4.} Further in Section \ref{gencan}, we study important cases, where Darboux matrices and fundamental solutions may be constructed
explicitly. The recovery of the generalised eigenfunction $\Pi(x)$, that is, solving of the matrix differential equation \eqref{c5}
is the main part of the explicit construction of the Darboux matrices and fundamental solutions.
We express this generalised eigenfunction in terms of the  square root of the matrix eigenvalue (see Proposition \ref{dem}).
In this way, a much wider class of generalised eigenfunctions and explicit solutions  (than in \cite{ALS17}) is constructed
even for $r=1$.

Matrix roots (especially, square roots of matrices) are actively studied and find various applications during more than one and half centuries
(see the paper \cite{Caley} from the year 1858). In particular,  structured roots of structured matrices and constructive
procedures to obtain these roots are of interest. Among important works on this topic one could mention \cite{Cross, Gant,  Wed}.
See also interesting recent works
\cite{Choi, Gaw} and references therein.
In the next section ``Preliminaries", we recall the construction of the matrix roots. Moreover, we modify
this construction so that the roots  possess certain commutation properties which are required in several
GBDT cases \cite{ALSa1, ALSstring, ALSJDE21}.

{\bf 5.} In Section \ref{Dyn}, we apply the results of Section \ref{gencan} in order to construct solutions
of the dynamical systems \eqref{I3} and study some of their properties.

Appendix \ref{GBDT} is dedicated to the general GBDT results from \cite{SaSaR} (see Theorem~\ref{GDGBDTrd})
and to Corollary \ref{CyGBDT} (for symmetric $S$-nodes), which we use in Section~\ref{gencan}.

Appendix \ref{dD} deals with another important application of matrix roots, that is, an application
to discrete Dirac systems and Verblunsky-type coefficients (see \cite{ALS2022} and references therein).

{\it Notations.} Many notations were explained above.  The notation 
$\sigma(A)$ stands for the spectrum of the matrix $A$,
$A^*$ stands for the conjugate transpose of  $A$, and the matrix inequality $S>0$ ($S\geq 0$)  means  that the eigenvalues of the matrix $S=S^*$ are positive
(nonnegative). In a similar  way, one interprets  the inequalities $S<0$ and $S\leq 0$.
The notation $\diag\{\cla_1, \cla_2, \ldots, \cla_s\}$ means diagonal matrix with the entries (or blocks)
$\cla_1, \cla_2, \ldots$ on the main diagonal.

%%%%%%%%%%%%%%%%%%%%%%%%%%%%%%%%%%%%%%%%%%%
\section{Preliminaries} \label{prel}
\setcounter{equation}{0}
{\bf 1.} An $\ell$-th root $Q$ $(\ell>1, \,\, \ell \in \BN)$ of  some invertible $n\times n$ matrix $A$ is of interest in many applications. (Such roots are constructed, e.g.,
in \S 6 of \cite[Chapter VIII]{Gant}). If one takes a root of $f(A)$ instead of $A$, it is natural to assume that $A$ and $Q$
still commute (see, e.g., \cite[Proposition 3.1]{ALSJDE21} where such commutation is required). In order to achieve  the equality $AQ=QA$ in this case, the construction in \cite[Chapter VIII, \S6]{Gant}
should be slightly modified. Namely, one should consider the matrices of the form 
\begin{align} &       \label{p1}
F=\mu I_p+T \quad (\mu\not=0, \quad T\in \cln\clt),
 \end{align} 
{\it where $\cln\clt$ is the class of upper triangular nilpotent Toeplitz matrices}, while Jordan cells are dealt with in \cite{Gant}.
Clearly, $p\times p$ matrix $T \in \cln\clt$ admits representation
\begin{align} &       \label{p1+}
T=\sum_{i=1}^{p-1}t_i\cls_i,
 \end{align} 
where $\cls_i$ are shift matrices:
\begin{align}& \label{p!}
 \cls_{i}:=\{\delta_{k-r+i}\}_{k,r=1}^p \quad (i>0), \quad \cls_i \cls_k =\cls_{i+k},
\end{align}
$\delta_s$ is Kronecker delta, and $\cls_{i}=0$ for $i \geq p$. It is known (and follows immediately from $\cls_i \cls_k =\cls_{i+k}$)
that the matrices from $\cln\clt$ commute.

Considering the branch of $h(\la)=(\la+\mu)^{1/\ell}$ (for small $\la$) near one of the roots $\mu^{1/\ell}$ and taking the partial sum of the first $p$ terms of the Taylor expansion of $(\la+\mu)^{1/\ell}$
at $\la =0$,
we set
\begin{align}       \nn
g(\la)=&\mu^{1/\ell}+\frac{1}{\ell }\mu^{\left(\frac{1}{\ell }-1\right)}\la+\frac{1}{\ell }\left(\frac{1}{\ell }-1\right)\mu^{\left(\frac{1}{\ell }-2\right)}\frac{\la^2}{2}+\ldots
\\ &  \label{p2}
+
\frac{1}{\ell }\left(\frac{1}{\ell }-1\right)\ldots \left(\frac{1}{\ell }-p+2\right)\mu^{\left(\frac{1}{\ell }-p+1\right)}\frac{\la^{p-1}}{(p-1)!}.
 \end{align} 
Then, we have $g(\la)^{\ell}=\mu +\la+O\big(\la^{p}\big)$ for $\la\to 0$. Thus, $g(\la)^{\ell}$ and $\mu +\la$ {\it coincide on the spectrum of
$T$} (in the terminology of \cite[Chapter V, \S1]{Gant}), and so 
\begin{align} &       \label{p3}
g(T)^{\ell}=\mu I_p+T=F .
\end{align}
Let the minimal polynomial of a matrix $A$ have the form $\prod_{i=1}^s(\la-\la_i)^{p_i}$.
Recall that $g(\la)$ is determined on the spectrum of this matrix $A$ (i.e., $g(A)$ exists)  if $g(\la_i )$ and its derivatives
up to $g^{(p_i)}(\la_i)$ exist in the points $\la_i$ \\ $(1\leq i\leq s)$ of the spectrum of $A$.

{\bf 2.} Given some $n\times n$ matrix $A$ ($\det A\not=0$), we consider its representation
\begin{align} &       \label{p4}
A=uA_{\clj}u^{-1}, \quad A_{\clj}=\diag\{\cla_1, \cla_2, \ldots, \cla_s\}, \quad \cla_i=\mu_i I +\cls_1,
\end{align}
where $A_{\clj}$ is a Jordan normal form of $A$,  $\cla_i$ are Jordan cells (of order $p_i$),
$I$ and $\cls_1$ have the same order $p_i$ as $\cla_i$, and  $I$ is the identity matrix. (We omit sometimes
the orders of $I$ and $\cls_k$ in the notations.) We note that $g(\la)$ depends, in fact, on two additional parameters:
$g(\la)=g(\la, \mu, k)$ $(0 \leq k<\ell)$. Here, $k$ determines the chosen root $\mu^{1/\ell}=\E^{2\pi \I k/\ell}\xi_0$
and $\xi_0$ is some arbitrarily fixed root.

The calculation of square roots $Q$ such that $AQ=QA$ and either $Q^2 = (A-cI_n)^2 - |a|^2 I_n$
or $Q^2 = (A-cI_n)^2 + |a|^2 I_n$ was required for the construction of Darboux matrices and explicit solutions
of Dirac systems in \cite[Proposition 3.1]{ALSJDE21}. Similar problems appear in \cite[Proposition 4.1]{ALSstring}.
\begin{Pn}\label{PnRoot} Let Jordan normal form of the $n\times n$ matrix $A$  $(\det A\not=0)$ 
be given by \eqref{p4} and let the function $f$ be defined
on the spectrum of $A$. 

Then, the matrices $Q$ of the form 
\begin{align} &       \label{p5}
Q=u\, \diag\{\clb_1, \clb_2, \ldots, \clb_s\}\, u^{-1}, \quad \clb_i=g\big(f(\cla_i)-f(\mu_i) I, f(\mu_i), k_i\big),
\end{align}
where $k_i$ may take fixed arbitrary  values between $0$ and $\ell-1$, satisfy relations
\begin{align} &       \label{p6}
Q^{\ell}=f(A), \quad AQ=QA.
\end{align}
\end{Pn}
\begin{proof}. It is easy to see that 
\begin{align} &       \label{p7}
f(A)=u\, \diag\{f(\cla_1), f(\cla_2), \ldots, f(\cla_s)\}\, u^{-1}.
\end{align}
Moreover, in view of the last equalities in \eqref{p!} and \eqref{p4} we have  
$$f(\cla_i)-f(\mu_i) I \, \in \, \cln\clt.$$
Hence, formulas \eqref{p3}, \eqref{p5} and \eqref{p7} yield the first equality in \eqref{p6}.

Finally, since $\clb_i-\big(f(\mu_i)\big)^{1/\ell} I\, \in \, \cln\clt$, the last relation in \eqref{p4}
implies that $\cla_i\clb_i=\clb_i\cla_i$ and the second equality in \eqref{p6} follows.
\end{proof}
\begin{Rk} \label{RkExp} Clearly, if $u$ and $A_{\clj}$ are given, the roots $Q$ of $f(A)$ are constructed in Proposition \ref{PnRoot}
explicitly.
\end{Rk}
\begin{Rk} \label{RkComR} It follows from $\clb_i-\big(f(\mu_i)\big)^{1/\ell} I\, \in \, \cln\clt$ that different roots of $f(A)$ constructed
in Proposition \ref{PnRoot} commute.
\end{Rk}
\begin{Rk} \label{RkNonCom} It is easy to see that the equality $AQ=QA$ does not necessarily hold for all the roots
of $Q^{\ell}=f(A)$ $($and so not necessarily all the roots are constructed in in Proposition \ref{PnRoot}$)$. For instance, setting 
$n=3$, $A=I_3+\cls_1$, $f(\la)=(\la-1)^2+4$, 
and $\ell=2$, we can define $u$ and one of the roots $Q$
as follows:
\begin{align} &       \label{p8}
u=u^{-1}=\begin{bmatrix} 1 & 0 & 0\\ 0&0& 1\\ 0 & 1 & 0\end{bmatrix}, \quad Q=\begin{bmatrix} 2 & 0 & 1/4 \\ 0&-2 & 0\\ 0 & 0 & 2\end{bmatrix} 
\quad \big(Q^2=(A-I_3)^2+4I_3\big).
\end{align}
One can easily check that $\cls_1Q\not=Q\cls_1$ and so $AQ\not=QA$.
\end{Rk}

{\bf 3.} Another (related to $AQ=QA$) commutation property is required for the explicit construction of solutions
of sigma models and gravitational equations in \cite{ALSa1}. Namely (after a change of notations for the ones used in this paper),
it is required in \cite[(B.1) and (B.2)]{ALSa1} that
\begin{align}& \label{p9}
Q(z)^{\ell}=A -z I_n, \quad 
Q(z_1)Q(z_2)=Q(z_2)Q(z_1) \quad (z_1,z_2\in \BR),
\end{align} 
where $z,z_1,z_2 \not\in \sigma(A)$ and $\ell=2$.

In order to solve \eqref{p9}, we  fix $u$ (independently  of $z$) in the procedures to find the roots $Q(z)$ of $A-z I_n$ in Proposition \ref{PnRoot}.
Then, we  find the blocks $\clb_i(z_1)$ setting $f(\la)=\la -z_1$ and find the blocks $\clb_i(z_2)$ setting $f(\la)=\la -z_2$.
The obtained blocks are upper triangular Toeplitz matrices. Thus, the second (commutation)
equality in \eqref{p9} holds. In this way,  the following corollary of  Proposition \ref{PnRoot} is obtained.
\begin{Cy} Let $\ell>1, \,\, \ell \in \BN$ and assume that  $z,z_1,z_2 \not\in \sigma(A)$. Then, 
the roots $Q(z)$ satisfying \eqref{p9} are constructed using procedures from Proposition~\ref{PnRoot}.
\end{Cy}

\section{GBDT for generalised Hamiltonian systems} \label{gencan}
\setcounter{equation}{0}
{\bf 1.} Recall that we consider generalised Hamiltonian systems of the form
\begin{align} &       \label{c1}
y^{\prime}(x,z)=G(x,z)y(x,z), \quad G(x,z)=\I \sum_{k=1}^{r} (z-c_k)^{-1}j\b_k(x)^*\b_k(x),
\end{align}
where \eqref{I2} holds. Systems \eqref{c1} are considered on some interval $\cli$ $(0\in \cli)$.

Given initial system of the form \eqref{c1}, we construct GBDT-transformed system, its fundamental solution
and the corresponding Darboux matrix using general results on GBDT
from Appendix \ref{GBDT}. 
For this purpose, we set $J=j$ in \eqref{A1} (and further in Appendix \ref{GBDT}) and fix some triple of matrices
$\{A, \, S(0),\, \Pi(0)\}$ such that
\begin{align} &       \label{c2}
A, S(0)\in \BC^{n\times n}, \quad \Pi(0)\in \BC^{n\times m}, \quad n\in \BN, \quad c_k\not\in\sigma(A)\,\, (1\leq k\leq r),
\\ &       \label{c3}
AS(0)-S(0)A^*=\I \Pi(0)j\Pi(0)^*, \quad S(0)=S(0)^*.
\end{align}
We note that  some notations here are different from Appendix \ref{GBDT}. In particular, the summation in \eqref{c1} (from $k=1$ till $r$)
is substituted by the summation from $s=1$ till $\ell$ in \eqref{GDg1}. Lower indices below are enumerated sometimes with commas in order to separate one index from another.
Clearly, the equalities on $q_{1,k}$ in \eqref{A3} hold for
 \begin{align} &       \label{c4}
q_{1,k}=-\I jH_k(x), \quad H_k(x):=\b_k(x)^*\b_k(x) .
\end{align}

According to \eqref{c1}--\eqref{c4}, the conditions of Theorem \ref{GDGBDTrd} and Corollary~\ref{CyGBDT}
for GBDT determined by the triple $\{A, \, S(0),\, \Pi(0)\}$ are fulfilled. Taking into account the form
of $G(x,z)$ in \eqref{c1},  equalities \eqref{A2} and the first equality in \eqref{A4}, we rewrite  formulas
\eqref{GDg3}, \eqref{GDg5'}, \eqref{GDg6}, and \eqref{GDJ5} (in this order) in the form
 \begin{align}   &       \label{c5}
 \Pi^{\prime}(x)=-\I\sum_{k=1}^{r} (A -c_k I_n)^{-1}\Pi(x)jH_k(x),
\\  &       \label{c6} 
S^{\prime}(x)=-\sum_{k=1}^{r} (A -c_k I_n)^{-1} \Pi (x)jH_k(x)j\Pi(x)^*(A^* -c_k I_n)^{-1},
\\  &       \label{c7} 
AS(x)-S(x)A^*=\I \Pi(x)j\Pi(x)^*,
\\ &       \label{c8}
w_A(x,z)=I_m-\I j\Pi(x)^*S(x)^{-1}(A-z I_n)^{-1}\Pi(x).
\end{align}
The matrix functions $\Pi(x)$ and $S(x)$ are uniquely determined by the triple $\{A, \, S(0),\, \Pi(0)\}$ and systems \eqref{c5}, \eqref{c6}.
\begin{Rk} \label{RkExpl} We will express the GBDT-transformed system via the transfer function $w_A$, which, in turn, is expressed via $\Pi(x)$
and $S(x)$. Moreover, when $\sigma(A)\cap\sigma(A^*)=\emptyset$, $S(x)$ is uniquely recovered from  \eqref{c7}.
The Darboux matrix and fundamental solution of the GBDT-transformed system are expressed via $w_A$ as well
$($see Appendix \ref{GBDT} and the text below in this section$)$. 
Summing up, we see that explicit construction of $\Pi(x)$ is one of the main components of this theory $($and construction of $S(x)$ and $w_A(x,z)$
easily follows$)$. Matrix roots are essential for the construction of $\Pi(x)$  as demonstrated in Proposition \ref{dem}.
\end{Rk}
Using \eqref{c4} and \eqref{GDg11}, we see that the transformed system (for the initial system \eqref{c1}) has the form
\begin{align} &       \label{c9}
\wt y^{\prime}(x,z)=\I \sum_{k=1}^{r} (z-c_k)^{-1}j\wt H_k(x)\wt y(x,z), 
\\  \nn &
\wt H_k=H_k+H_kY_{k,-1}-jX_{k,-1}jH_k-jX_{k,-1}jH_kY_{k,-1}
\\ &       \label{c10} \quad \,
=(I_m-jX_{k,-1}j)H_k(I_m+Y_{k,-1}).
\end{align}
Relations \eqref{GDg13} and \eqref{GDJ5} (together with the formula \cite[(1.76)]{SaSaR} show that
 \begin{align} &       \label{c11}
I_m-X_{k,-1}(x)=w_A(x,c_k), \quad \big(I_m-X_{k,-1}(x)\big)\big(I_m+Y_{k,-1}(x)\big)=I_m.
\end{align}
On the other hand, formula \cite[(1.87)]{SaSaR} for symmetric $S$-nodes yields
 \begin{align} &       \label{c12}
w_A(x,c_k)^{-1}=jw_A(x,c_k)^*j.
\end{align}
Equalities \eqref{c4} and \eqref{c10}--\eqref{c12} imply that
 \begin{align} &       \label{c13}
\wt H_k(x)=\wt\b_k(x)^*\wt \b_k(x), \quad \wt \b_k(x):=\b_k(x) jw_A(x,c_k)^*j.
\end{align}
The matrices $\wt H_k(x)$ are nonnegative and may be considered as the transformed generalised Hamiltonians. 
Our next statement follows from Theorem~\ref{GDGBDTrd}.
\begin{Tm} \label{GBDTgc} Let the initial generalised Hamiltonian system   \eqref{c1} and 
the triple $\{A, \, S(0),\, \Pi(0)\}$ be given, where \eqref{I2}, \eqref{c2} and \eqref{c3} hold, 
and  $c_k\not=c_l$ for $k\not=l$. Then, the GBDT-transformed generalised Hamiltonian system 
is given by the formulas \eqref{c9} and \eqref{c13}, and $w_A(x,z)$ determined by \eqref{c5}, \eqref{c6} and \eqref{c8}
is the corresponding  Darboux matrix. If $S(x)$ is invertible on the interval $\cli$  where our systems
are considered, and $w(x,z)$ is a fundamental solution of the initial system, then the fundamental solution $\wt w(x,z)$ of the
transformed system has the form: 
\begin{align} &       \label{c14}
\wt w(x, z)=
w_{A}(x, z )w(x, z ).
 \end{align} 
 \end{Tm} 
 \begin{Rk}\label{RkS} According to \eqref{c6}, $S^{\prime}(x)\leq 0$. Thus, $S(x)<0$ in the case 
 $$S(0)<0,  \quad \cli=[0,a] \quad (a>0),$$
 and $S(x)>0$ in the case 
 $$S(0)>0, \quad \cli=[-a,0] \quad (a>0).$$
In both cases, $S(x)$ is invertible. 
 \end{Rk}
 {\bf 2.}  Let us consider some examples, where the conditions of  Theorem \ref{GBDTgc} hold.
 \begin{Ee}\label{Ee1}  Consider a simple example of the trivial initial Hamiltonians, that is,
 the case of $H_k(x)\equiv I_m$ $(1\leq k\leq r)$. Assuming for simplicity that $m_1>0$ and $m_2>0$, the case may be treated  similar $($in many details$)$ to the
 example with $r=1$ in \cite{ALS17}. 
 
 We partition $\Pi(0)$ into the $n\times m_1$ and $n\times m_2$ blocks $\vt_1$ and $\vt_2$, respectively: 
 $\Pi(0)=\begin{bmatrix} \vt_1 & \vt_2\end{bmatrix}$. Hence, it  follows from \eqref{c5} that 
 \begin{align} &       \label{E1}
\Pi(x)=\begin{bmatrix}\exp\Big\{-\I x\sum_{k=1}^{r} (A -c_k I_n)^{-1}\Big\}\vt_1 & \exp\Big\{\I x\sum_{k=1}^{r} (A -c_k I_n)^{-1}\Big\}\vt_2\end{bmatrix}.
 \end{align} 
 If $\s(A)\cap\s(A^*)=\emptyset$, there are unique matrices $C_i$ $(i=1,2)$ such that
 \begin{align} &       \label{E2} 
  AC_i-C_iA^*=\I \vt_i\vt_i^*.
  \end{align} 
  Moreover, in view of  \eqref{E1}, \eqref{E2} and equality $\s(A)\cap\s(A^*)=\emptyset$,  there is a unique solution
  $S(x)$ of \eqref{c7} which is  given by
   \begin{align}  \nn     
  S(x)=&\exp\Big\{-\I x\sum_{k=1}^{r} (A -c_k I_n)^{-1}\Big\}C_1\exp\Big\{\I x\sum_{k=1}^{r} (A^* -c_k I_n)^{-1}\Big\}
 \\ &  \label{E3} 
  -\exp\Big\{\I x\sum_{k=1}^{r} (A -c_k I_n)^{-1}\Big\}C_2\exp\Big\{-\I x\sum_{k=1}^{r} (A^* -c_k I_n)^{-1}\Big\}.
  \end{align}
  Here, equality \eqref{c7} for $S(x)$ of the form \eqref{E3} is checked by direct substitution.  In particular,
  we have $S(0)=C_1-C_2$.
  \end{Ee}
  
{\bf 3.} 
Next, we study a nontrivial example, where the matrix roots are involved.
 The following proposition is valid.
\begin{Pn}\label{dem}  Let the conditions of Theorem \ref{GBDTgc} hold. Assume that \\ $r=2$,
$m_1=m_2=p>0$, $\b_k(x)$ are constant matrices $\b_k$ $(k=1,2)$ and
\begin{align} &       \label{c15}
\b_kj\b_k^*=0 \quad (k=1,2), \quad \b_1j\b_2^*=I_p.
 \end{align} 
Then, the solution $\Pi(x)$ of  \eqref{c5} is given explicitly by the formulas
\begin{align} &       \label{c16}
\Pi(x)j\begin{bmatrix} \b_1^* & \b_2^*\end{bmatrix}=\begin{bmatrix} \Phi_1(x) & \Phi_2(x)\end{bmatrix};
\\ &       \label{c17}
\Phi_1(x)=\E^{\I x Q}h_1+\E^{-\I x Q}h_2, \quad h_1,\, h_2, \, \Phi_1(x) \in \BC^{n\times p},
\\ &       \label{c18}
 \Phi_2(x)=-(A-c_2I_n)Q\big(\E^{\I x Q}h_1-\E^{-\I x Q}h_2\big), \quad \Phi_2(x)\in \BC^{n\times p};
\\ &       \label{c19}
h_1+h_2=\Pi(0)j\b_1^*, \quad (A-c_2I_n)Q(h_2-h_1)=\Pi(0)j\b_2^*;
\\ &       \label{c20}
Q^2=(A-c_1I_n)^{-1}(A-c_2 I_n)^{-1}.
 \end{align}
\end{Pn}
\begin{proof}.  First note that 
$$\begin{bmatrix} \b_1 \\ \b_2\end{bmatrix}j\begin{bmatrix} \b_1^* & \b_2^*\end{bmatrix}=\begin{bmatrix} 0 & I_p \\ I_p & 0\end{bmatrix}.$$
Thus, $\begin{bmatrix} \b_1^* & \b_2^*\end{bmatrix}$ is invertible (and $\begin{bmatrix} \b_1^* & \b_2^*\end{bmatrix}^{-1}=\begin{bmatrix} \b_2 \\ \b_1\end{bmatrix}j$).
Hence, in view of \eqref{c15} and \eqref{c16}, formula \eqref{c5} is equivalent to the equalities
\begin{align} &       \label{c21}
\Phi_1^{\prime}(x)=-\I (A-c_2 I_n)^{-1}\Phi_2(x), \quad \Phi_2^{\prime}(x)=-\I (A-c_1 I_n)^{-1}\Phi_1(x).
\end{align} 
The first equality in \eqref{c21} is immediate from \eqref{c17} and \eqref{c18}. The second equality in
\eqref{c21} follows from \eqref{c17}, \eqref{c18}, and \eqref{c20}.  Now, we see that
$$\Pi(x)=\begin{bmatrix} \Phi_1(x) & \Phi_2(x)\end{bmatrix}\begin{bmatrix} \b_2 \\ \b_1\end{bmatrix}$$
satisfies \eqref{c5}.

Relations \eqref{c19} follow from \eqref{c16} and from the equalities \eqref{c17} and \eqref{c18} at $x=0$.
If $\Pi(0)$ is given, $h_1$ and $h_2$ are uniquely recovered from \eqref{c19}. (Vice versa, starting from
$h_1$ and $h_2$ we easily obtain $\Pi(0)$.)
 \end{proof}

\begin{Rk}\label{RkB} We note that commutation properties are not required from the square
root $Q$ in \eqref{c20}. The requirement $AQ=QA$ is convenient but not necessary in \cite{ALSstring, ALSJDE21}.
However, in using GBDT for the construction of explicit solutions of integrable nonlinear PDEs the commutation requirements
are often unavoidable $($see, e.g., paragraph {\bf 3} of Section \ref{prel}$)$.
\end{Rk}

\section{The case of several variables}\label{Dyn}
\setcounter{equation}{0}
{\bf 1.} Let the conditions of Theorem \ref{GBDTgc} be fulfilled. That is, let the corresponding set of
matrix functions $\b_k(x)$ and the triple $\{A, \, S(0),\, \Pi(0)\}$ be given.
It follows from Corollary  \ref{CyGBDT} that, for our case of symmetric $S$-nodes and system \eqref{c1}, 
formula \cite[(7.61)]{SaSaR} takes the form
\begin{align} &       \label{d1}
\Big(j\Pi(x)^*S(x)^{-1}\Big)^{\prime}=\I\sum_{k=1}^rj\wt H_k(x) j\Pi(x)^*S(x)^{-1}(A-c_k I_n)^{-1}.
 \end{align} 
 Relation \eqref{d1} yields our next theorem.
 \begin{Tm} \label{TmSV} Let a set of numbers $c_k$ and matrices $\b_k(x)$ $(1\leq k \leq r)$ as well as a triple
$\{A, \, S(0),\, \Pi(0)\}$ be given. Assume that relations \eqref{I2}, \eqref{c2} and \eqref{c3} hold
and that $\Pi(x)$ and $S(x)$ satisfy \eqref{c5} and \eqref{c6}, respectively. 

Then, in the points of invertibility of $S(x)$, the matrix function
 \begin{align} &       \label{d2}
\wt \psi(x,\ze_1,\ldots,\ze_r)=j\Pi(x)^*S(x)^{-1}\exp\Big\{\sum_{k=1}^r\ze_k(A-c_k I_n)^{-1}\Big\}
 \end{align} 
 satisfies the system
 \begin{align} &       \label{d3}
\frac{\p \wt \psi}{\p x}=\I\sum_{k=1}^r j \wt H_k(x)\frac{\p \wt \psi}{\p \ze_k},
\end{align}
where $\wt H_k$ is given by \eqref{c13}.
 \end{Tm}
\begin{Rk} Note that $\wt \psi$ is an $m\times n$ matrix function. Thus,
\eqref{d3} means that each column of $\wt \psi$ satisfies our transformed
dynamical system $($of several variables$)$.
\end{Rk}
\begin{Rk}\label{RkSi}
Under conditions of Theorem \ref{TmSV}, equalities \eqref{c12} and \eqref{c13} $($where $w_A$ is expressed via $\Pi$ and $S$ in
\eqref{c8}$)$ are fulfilled. Hence, one may rewrite \eqref{c13} in the form:
\begin{align} &       \label{d4}
j\wt H_k(x)=w_A(x,c_k)jH_k(x)w_A(x,c_k)^{-1},
 \end{align}  
that is,  $j\wt H_k(x)$ is linear similar to $j H_k(x)$.
\end{Rk}
{\bf 2.} The following conservation law is valid for the GBDT-transformations in this paper.
\begin{Pn}\label{PnCL} Let the conditions of Theorem \ref{TmSV} hold. Then,
\begin{align} &       \label{d5}
\big(\Pi(x)^*S(x)^{-1}\Pi(x)\big)^{\prime}=\sum_{k=1}^r\big(\wt H_k(x)-H_k(x)\big).
 \end{align}  
\end{Pn}
\begin{proof}.
Relations \eqref{c8} and \eqref{d1} yield
\begin{align} &       \label{d6}
\big(\Pi(x)^*S(x)^{-1}\big)^{\prime}\Pi(x)=\sum_{k=1}^r\wt H_k(x)\big(I_m-w_A(x,c_k)\big).
 \end{align}  
From \eqref{c5} and \eqref{d6},  it follows that
\begin{align}        \nn
\big(\Pi(x)^*S(x)^{-1}\Pi(x)\big)^{\prime}=&\sum_{k=1}^r\wt H_k(x)\big(I_m-w_A(x,c_k)\big)
\\ &
\label{d7}
+j\sum_{k=1}^r\big(w_A(x,c_k)-I_m\big)jH_k(x).
 \end{align}   
 Finally, the equalities \eqref{d4} and \eqref{d7} imply \eqref{d5}.
\end{proof}.
\begin{Rk}\label{RkDyn}
If $\Pi(x)$ is known explicitly $($as in Example \ref{Ee1} or in Proposition \ref{dem}$)$, we obtain explicit expressions for
$\wt \psi$ and system \eqref{d3}.
\end{Rk}
\appendix
\section{GBDT - general case} \label{GBDT}
\setcounter{equation}{0}
{\bf 1.} Following \cite[pp. 219--221]{SaSaR} (see also the references therein), let us  consider a general case of the first order 
system  depending rationally on the spectral parameter $z$, that is, the case of multiple poles with respect to $z$:
\begin{equation}       \label{GDg1}
y^{\prime}=Gy, \quad G(x,z)=-\Big(\sum_{k=0}^r z^k q_k(x)+\sum_{s=1}^{\ell} \sum_{k=1}^{r_s}(z-c_s)^{-k}
q_{sk}(x)\Big).
\end{equation}
Similar to \cite{SaSaR}, we assume that $x \in {\cal I}$, where ${\cal I }$ is an interval such that $0 \in {\cal I }$, and that the coefficients
$q_k(x)$ and $q_{sk}(x)$ are $m \times m$ locally summable matrix functions.
The GBDT of the initial system \eqref{GDg1} is  determined by a number $n \in \BN$,
by  $n \times n$ matrices $A_k$ ($k=1,2$) and $S(0)$, and by $n \times m$ matrices
$\Pi_k(0)$ ($k=1,2$). It is required that these matrices form an $S$-node, that is,
the matrix identity
\begin{equation}       \label{GDg2}
A_1S(0)-S(0)A_2=\Pi_1(0)\Pi_2(0)^*
\end{equation}
holds.  The  matrix functions $\Pi_k(x)$ are introduced via their values at $x=0$ and 
equations
\begin{eqnarray}       \label{GDg3}
&&\big(\Pi_1 \big)^{\prime}=\sum_{k=0}^r A_1^k \Pi_1 q_k+\sum_{s=1}^{\ell} \sum_{k=1}^{r_s}(A_1 -c_s I_n)^{-k}
\Pi_1 q_{sk},\\
      \label{GDg4}
&&\big(\Pi_2^* \big)^{\prime}=-\Big(\sum_{k=0}^r  q_k\Pi_2^*A_2^k+\sum_{s=1}^{\ell} \sum_{k=1}^{r_s}q_{sk}\Pi_2^*(A_2 -c_s I_n)^{-k}
\Big).
\end{eqnarray}
Compare  (\ref{GDg1}) with (\ref{GDg4}) to see that $\Pi_2^*$ can be viewed
as a generalized eigenfunction of the system $u^{\prime}=Gu$.

Matrix function $S(x)$ is introduced via $S(0)$ and via $S^{\prime}(x)$ given by the equality
\begin{eqnarray}       \nn S^{\prime}&=&\sum_{k=1}^r\sum_{j=1}^k A_1^{k-j} \Pi_1 q_k\Pi_2^*A_2^{j-1}-\sum_{s=1}^{\ell} \sum_{k=1}^{r_s}
\sum_{j=1}^k
(A_1 -c_s I_n)^{j-k-1} \\ \label{GDg5}
 && \times
\Pi_1 q_{sk}\Pi_2^*(A_2 -c_s I_n)^{-j}.
\end{eqnarray}
%%%%%%%%%%%%%%%%%%%%%%%%%%%%%%%%%%%%%%%%%%%%%%%
Equations   \eqref{GDg2}--\eqref{GDg5})  yield the
identity
\begin{equation}       \label{GDg6}
A_1S(x)-S(x)A_2=\Pi_1(x)\Pi_2(x)^*, \quad x\in {\cal I}
\end{equation}
of which \eqref{GDg2} is a particular case of $x=0$. For initial system \eqref{GDg1}, the GBDT-transformed system
is
\begin{equation}       \label{GDg9}
\wt y^{\prime}=\wt G \wt y, \quad \wt G(x,z)=-\Big(\sum_{k=0}^r z^k \wt q_k(x)+\sum_{s=1}^{\ell} \sum_{k=1}^{r_s}(z-c_s)^{-k}
\wt q_{sk}(x)\Big),
\end{equation}
where the transformed coefficients $\wt q_k$ and $\wt q_{sk}$ are given by the formulas
\begin{align}       \label{GDg10}&
\wt q_k=q_k-\sum_{j=k+1}^r\Big(q_jY_{j-k-1}-X_{j-k-1}q_j+\sum_{i=k+2}^jX_{j-i}q_j
Y_{i-k-2}\Big), 
\\
  \label{GDg11} &
\wt q_{sk}=q_{sk}
+\sum_{j=k}^{r_s}\Big(q_{sj}Y_{s,k-j-1}-X_{s,k-j-1}q_{sj}-\sum_{i=k}^jX_{s,i-j-1}q_{sj}
Y_{s,k-i-1}\Big),
\end{align}
and the matrix functions $X_k(x)$, $Y_k(x)$  and $X_{sk}(x)$, $Y_{sk}(x)$ 
have the form
\begin{align}    & \label{GDg12}
 X_k=\Pi_2^*S^{-1}A_1^k\Pi_1, 
\quad Y_k=\Pi_2^*A_2^kS^{-1}\Pi_1 ,
\\ &  \label{GDg13}
X_{sk}=\Pi_2^*S^{-1}(A_1-c_s I_n)^k\Pi_1, \quad Y_{sk}=\Pi_2^*(A_2-c_s I_n)^kS^{-1}\Pi_1.
\end{align}
Our theorem below shows that the transfer function
\begin{align} &       \label{GDJ5}
w_A(z)=I_m-\Pi_2^*S^{-1}(A_1-z I_n)^{-1}\Pi_1,
 \end{align} 
(that is, the transfer matrix function in Lev Sakhnovich form \cite{SaL1, SaL2})
is the Darboux matrix for the systems \eqref{GDg1} and \eqref{GDg9}.
\begin{Tm} \label{GDGBDTrd} \cite{SaSaR} Let the first order system   \eqref{GDg1} and five
matrices $S(0)$, $A_1$, $A_2$ and $\Pi_1(0)$, $\Pi_2(0)$   be given. Assume that
the identity  \eqref{GDg2} holds and that $\{c_s\}\cap \s(A_k)=\emptyset$ $(k=1,2)$. Then, in the
points of invertibility of $S$, the transfer matrix function $w_A(x,z)$ given by  
\eqref{GDJ5}, where $S(x)$ and $\Pi_k(x)$ are determined by \eqref{GDg3}--\eqref{GDg5},
 satisfies the equation 
\begin{align} &       \label{A0}
w_{A}^{\prime}(x, z )=
\widetilde{G}(x, z ) w_{A}(x, z )- w_{A}(x, z )
G(x, z),
 \end{align} 
where $\wt G$ is determined
 by the formulas  \eqref{GDg9}--\eqref{GDg13}.
 \end{Tm} 
Compare system \eqref{A0} with the system \eqref{I0} in introduction.

We note that the  transfer matrix function $w_A(z)$ given by \eqref{GDJ5}  takes roots in M.  Liv\v{s}ic's
characteristic matrix function (see  the discussion in \cite{SaL21}).

{\bf 2.} In this paper, we are interested in the case of the {\it symmetric $S$-nodes}, that is, in the case
\begin{align} &       \label{A1}
 A_2=A_1^*, \quad S(0)=S(0)^*,  \quad \Pi_2(0)=\I J\Pi_1(0)^*, \quad J=J^*=J^{-1}.
 \end{align} 
In this case, we also change notations and set
\begin{align} &       \label{A2}
 A_1=A_2^*=A,   \quad \Pi_1(x)=\Pi(x).
 \end{align}  
 Moreover, we assume that
 \begin{align} &       \label{A3}
Jq_k^*J=-q_k, \quad Jq_{sk}^*J=-q_{sk}, \quad c_s\in \BR \quad (1\leq s\leq \ell, \quad 1\leq k\leq r_s).
 \end{align}  
The following corollary is valid under conditions above.
\begin{Cy} \label{CyGBDT} Let the conditions of Theorem \ref{GDGBDTrd} and relations \eqref{A1}--\eqref{A3} hold.
Then, we have
\begin{align} &       \label{A4}
\Pi_2(x)^*=\I J \Pi(x)^*, \quad S(x)=S(x)^*.
 \end{align} 
\end{Cy}
\begin{proof}.  It is immediate from \eqref{GDg3}, from \eqref{A1}, and from $c_s\in \BR$ that
\begin{align}       \label{GDg3'}
\big(\I J \Pi_1^* \big)^{\prime}
=\sum_{k=0}^r Jq_k^*J \big(\I J \Pi_1^*\big) A_2^k+\sum_{s=1}^{\ell} \sum_{k=1}^{r_s}
Jq_{sk}^*J  \big(\I J \Pi_1^*\big) (A_2 -c_s I_n)^{-k}.
\end{align}
Hence, in view of \eqref{GDg4}, \eqref{A3} and the third equality in \eqref{A1}, $\I j \Pi_1(x)^*$ and $\Pi_2(x)^*$ satisfy the same
first order linear differential system and coincide at $x=0$. Thus, we obtain the first equality in \eqref{A4}.
Therefore, taking into account \eqref{A2}, we rewrite \eqref{GDg5} in the form
\begin{align}       \nn S^{\prime}=&\sum_{k=1}^r\sum_{j=1}^k A^{k-j} \Pi (\I q_kJ)\Pi^*(A^*)^{j-1}-\sum_{s=1}^{\ell} \sum_{k=1}^{r_s}
\sum_{j=1}^k
(A -c_s I_n)^{j-k-1} \\ \label{GDg5'}
 & \times
\Pi (\I q_{sk}J)\Pi^*(A^* -c_s I_n)^{-j}.
\end{align}
Changing the order of summation, that is,  setting 
$$k-j=p-1, \quad Q_k=\Pi (\I q_kJ)\Pi^*,  \quad Q_{sk}=\Pi (\I q_{sk}J)\Pi^*,$$ 
we have
\begin{align} \nn       
\sum_{j=1}^k A^{k-j} Q_k(A^*)^{j-1}&=\sum_{p=1}^k A^{p-1}Q_k(A^*)^{k-p}
\\ & \label{A5}
=\frac{1}{2}\sum_{j=1}^k \big(A^{k-j} Q_k(A^*)^{j-1}+A^{j-1}Q_k(A^*)^{k-j}\big),
\end{align}  
and in, in a similar way,
 \begin{align}
 \nn &
\sum_{j=1}^k
(A -c_s I_n)^{j-k-1}
Q_{sk}(A^* -c_s I_n)^{-j}
\\ & \nn
=\frac{1}{2}\sum_{j=1}^k\big((A -c_s I_n)^{j-k-1}
Q_{sk}(A^* -c_s I_n)^{-j}
\\ & \qquad \label{A6}
+(A -c_s I_n)^{-j}
Q_{sk}(A^* -c_s I_n)^{j-k-1}\big).
\end{align}  
 Using \eqref{A3}, it is easy to see that $Q_k=Q_k^*$, $Q_{sk}=Q_{sk}^*$, and so the right-hand sides of \eqref{A5}
 and \eqref{A6} are Hermitian. It follows that the right-hand side of  \eqref{GDg5'} is Hermitian as well, that is,
 $S^{\prime}$ is Hermitian. Since $S(0)=S(0)^*$ and $S^{\prime}$ is Hermitian, the second equality in \eqref{A4}
 is also valid.
\end{proof}
In the case of symmetric $S$-nodes, formula \eqref{A1} implies that GBDT is determined by the triple $\{A, S(0), \Pi(0)\}$ of matrices.
\section{Matrix roots and discrete Dirac systems} \label{dD}
\setcounter{equation}{0}
Discrete Dirac systems
\begin{align} &\label{B1}
y_{k+1}(z)=(I_m-\frac{ \I}{ z}
j C_k)
y_k(z), \quad C_k>0, \quad C_k j C_k=j,
\end{align} 
present (in the case of $m_1=m_2$) a fruitful approach to Toeplitz matrices and measure theory (see \cite{ALS2022} and references therein), which completes in an
interesting way the classical approach via Szeg\H{o} recurrences. (Recall that $j$ is introduced in \eqref{I2}.)

The positive roots of the matrices $C_k$ play an essential role in the study of
discrete Dirac systems. Using \cite[Propositions 2.1, 2.4]{FKRS} (and completing the proof of \cite[Proposition 2.1]{FKRS} with the fact that there is a unique positive root $D^{1/\ell}>0$
of the matrix $D>0$, see \cite[(62), p. 232]{Gant}), we have the following proposition.
\begin{Pn}\label{PnD} Let some $m\times m$ matrix $C>0$ satisfy condition $CjC=j$. Then, there is a unique positive root $C^{1/\ell}$ $(\ell\in \BN)$
and this root has the property $C^{1/\ell}jC^{1/\ell}=j$. The matrix $C$ $($and $C^{1/\ell}$ as well$)$ admits so called Halmos extention representation: 
\begin{align}\nn &
C= {\mathrm{diag}}\Big\{
\big(
I_{m_1}-  \rho  \rho^* \big)^{-\frac{1}{2}}, \, \,
\big(I_{m_2}- \rho^* \rho\big)^{-\frac{1}{2}}\Big\}\left[
\begin{array}{cc}
I_{m_1} &  \rho \\   \rho^* & I_{m_2}
\end{array}
\right], \quad \rho^*\rho <I_{m_2},
\end{align}
where $\rho$ is an $m_1\times m_2$ matrix.
\end{Pn} 
The properties of Halmos extensions are well-known (see, e.g., \cite{Dy0, DFK}).

{\bf Acknowledgments}  {This research    was supported by the
Austrian Science Fund (FWF) under Grant  No. Y-963.}

\end{document}